# Faith Believes, Hope Expects: The Impact of Calvin's Theology on the Mathematics of Chance

*Timothy C. Johnson*[1]

## Abstract

This paper attributes the sudden emergence of mathematical probability and statistics in the second half of the seventeenth century to Calvin's Reformed theology. Calvin accommodated Epicurean chance with Stoic determinism and synthesised *phronesis/prudentia*, founded personal experience and employed to deal with *tyche/fortuna*, and *episteme/scientia*, universal knowledge. This meant that matters of chance, which had previously been considered too particular for mathematical treatment, became part of *episteme/scientia*. Clear evidence of the significance of Calvin in mathematics is in the facts that Huygens considered using the word 'hope' to describe mathematical expectation and French mathematics still uses *espérance* for mathematical expectation. Calvin asserted that Hope represented a universal, objective and indubitable idea making it characteristic of mathematics. The argument is built on a review of how the ideas of Hope, Faith and Prudence have evolved in European thought that highlights Calvin's innovations. The conclusion identifies contemporary issues in the application of mathematics in society that are illuminated in light of Calvin's doctrine.

## Introduction

Hacking (2006) stimulated a revival of interest in the history of mathematical approaches to uncertainty by investigating why the concept of mathematical probability did not exist until the mid-seventeenth century. The answer he provided was that before the Renaissance probability related to matters of opinion, with opinions being founded on authority. During the Renaissance, the idea of non-deductive inference, 'internal evidence' from 'natural signs' that supported a conclusion without necessitating it, emerged. A criticism of Hacking's account is that there is little evidence to support the assertion that the idea of 'evidence' was not present before the Renaissance (van Brankel, 1976) (Howson, 1978) (Garber & Zabell,

[1] Department of Actuarial Mathematics and Statistics and the Maxwell Institute for Mathematical Sciences, Heriot-Watt University, Edinburgh, EH14 4AS, UK. T.C.Johnson@hw.ac.uk

1979), (Brown, 1987) (Franklin, 2001, p. 373). This does not undermine Hacking's general thesis that there was a shift in the interpretation of signs that changed people's understanding of 'evidence' that led to mathematical probability. It does mean that an explanation still needs to be offered as to why there was a dramatic shift in the seventeenth century that resulted in a mathematical approach to chance. This paper provides an explanation.

Hacking had a clear focus on mathematical approaches to chance, framing his argument in the context of Todhunter (2014). Critics of Hacking's thesis have a broader view of evidence and probability before 1650, for example Franklin (2001) presents a "history of rational methods of dealing with uncertainty". Like Hacking, this paper focuses on mathematical approaches to uncertainty and defines an approach as mathematical if it has specific properties (d'Alembert, 2009) (Heims, 1980, p. 83). Firstly, it should be based on abstractions or generalisations. Mathematical models might be applied to particular circumstances, but a mathematical approach should be founded on generalised ideas. A resultant property is that a mathematical approach to uncertainty should be universally applicable. An implication of universality is that a mathematical result should be 'objective' in the sense that it is not specific to a particular, 'subjective', viewpoint. Finally, a mathematical approach to uncertainty should seek to establish something indubitable about uncertainty. The abstraction to a universal idea needs to result in a framework that is believed to be true. These imply that a quantitative method for solving a problem is not 'mathematical' simply because it involves numbers or equations, mathematisation requires universality and confidence.

On this basis, we shall argue that Calvin's Reformed theology, developed in the mid-sixteenth century, created the conceptual framework for the mathematical treatment of chance. Calvin regarded the mundane experience of chance as real and rejected attempts by both rationalism and mysticism to control chance, employing an Epicurean argument that God's will was indiscernible, to maintain divine omnipotence. The Epicurean aspect of Reformed theology was eclipsed by a resurgence of Stoic determinism with the ascendency of Providentialism (Puritanism) in the early seventeenth century. While Hacking identifies Paracelsus (1493-1541) as important in developing the idea of 'internal evidence' (Hacking, 2006, pp. 35-45) we shall argue that 'Experimental Calvinism', a feature or Puritanism, was a more substantial basis for a revolution in the idea of evidence from signs, because it was more widespread. While the influence of Epicurean ideas on Reformed doctrine waned, Calvin's theology had a persistent impact removing the distinctions between different types of knowledge, notably

*phronesis*/*prudentia* and *episteme*/*scientia*. The synthesis of the two meant that judgement clouded by uncertainty could be addressed using universal scientific laws and, in particular, mathematics.

The argument in the paper was not developed from the perspective of Reformed theology. Rather, the original motivation for this paper is in Hacking's (2006, pp. 93-95) observation that Huygens had considered using the Latin word for 'hope', *spes*, for 'expectation' in the mathematical sense, and the French use of *espérance* ('hope') in mathematics is widespread today. To answer why hope played a part in mathematical expectation, the evolution of hope, from Pre-Socratic times to the Renaissance, was traced. This highlighted that both hope and faith were originally connected to the mitigation of uncertainty. With the emergence of philosophy, prudence, or practical wisdom, became the most important intellectual tool for addressing chance. Prudence was considered separate from scientific knowledge and this distinction remained in pre-Reformation Christian thinking. The significance of Calvin to the mathematisation of uncertainty became apparent after exploring the basis of Calvin's formula "faith believes, hope expects" and this revealed how Calvin transformed the notion of Hope from the "uncertain certainty" of the Scholastics to an absolute certainty. While the starting point had been the relevance of hope, it became apparent that Calvin's effect on the concept of prudence was more significant in the long-term.

The novelty of the paper is in its engagement with Reformed theology of the seventeenth century. Scholarship regarding Calvin usually focuses on theology rather than the secular and when the relationships between mathematics and Calvinist doctrine have been investigated, they focus on the impact of the mathematical theory on religious doctrine, such as Reeves (2015). When the role of religion on the development of probability and statistics has been considered, it is seen as inhibiting. Kendall (1956, pp. 11-12) argued that pagan Greeks and Romans recognised the significance of chance, but this awareness disappeared with the advent of Christianity, which saw the "finger of God in everything". Hacking observed that

> Europe began to understand concepts of randomness, probability, chance and expectation precisely at that point in its history when theological views of divine foreknowledge were being reinforced (Hacking, 2006, pp. 2-3)

However, Hacking, whose aim was to understand the emergence of probability from a perspective of social constructivism (Hacking, 2006, pp. xiii-xvii), focussed on material medical practices and he overlooked the significance of Epicureanism in Calvin's theology, regarding religion from the more Stoic perspective of Providentialism. Franklin pays more attention to theology but, with a commitment to Catholicism (Franklin, 2006, pp. 65-81) (2009), does not engage with Reformed doctrine. Daston finds little evidence that theology contributed to the development of probability (Daston, 1998, pp. 38, 130-137, 151-156). While she highlights the significance of life insurance in the establishment of mathematical approaches to chance, she ignores the fact that the first mathematically managed insurance scheme, the Scottish Ministers' Widows' Fund, was created by Reformed clergymen.

Daston (1980, pp. 241-242) identifies a distinction between a 'moral' approach to probability, centred on equity and justice, which had dominated thinking before 1600, and a 'prudential' approach, which "weighed individual possibilities for profit or loss with an eye to securing an advantage", adapting Hacking's epistemic-aleatory distinction into a moral-prudential distinction. The 'prudential' approach to probability replaced the 'moral' approach with "the advent of a new model of explanation for the social sciences which emphasized social regularities rather than individual rationality, coupled with recognition of the independence of mathematical probability from its applications" (Daston, 1980, p. 235). In Daston (1998) the argument becomes one of a transformation of the subjectively reasoning *l'homme éclair* into the objectively calculating *l'homme moyen*. This paper investigates the move from a 'moral' to 'prudential' approach to probability and statistics in terms of the transformation of Prudence from the principal Classical virtue into a "stodgy", amoral, concept (den Uyl, 1991), (Hariman, 2003). This paper argues that Calvin's synthesis of *phronesis*/*prudentia* and *episteme*/*scientia* motivated the diminution of Prudence that drove the transformation that Daston describes.

One obvious fact connecting Reformed theology to the development of probability and statistics is that almost all the leading figures involved shared, or were raised in, the same theological perspective of Augustine's doctrine on predestination and salvation by faith alone, the basis of Reformed doctrine: Pascal, the Bernoullis, de Moivre, Montmort, Cramer, Bayes, Arbuthnot, Graunt – with only Pascal and Montmort not being Reformed. The only individuals who did not fit the pattern were Cardano and Fermat, and Fermat had strong connections to French Calvinism. Contemporary mathematicians who had no close links to Augustinian doctrine, such as Gassendi, Descartes, Hobbes, Newton or Leibnitz, made little or no contribution to mathematical probability.

The argument begins by summarising pre-Christian approaches to faith, hope and prudence. It then gives an account of Christian perspectives, describing Scholastic attitudes. This aims to provide clarity when considering the relevance of Scholastic ideas to probability, central to Franklin's case. The section then gives an account of Calvin's theology, highlighting how Calvin was required to innovate and a description of later Puritan Providentialism. The section finishes by summarising how attitudes to prudence changed from Machiavelli, a contemporary of Calvin and who regarded prudence in the Classical manner, to Hobbes and Spinoza, who, strongly influenced by Reformed religion, sought more scientific approaches to uncertainty. The purpose of this section is to provide sufficient evidence for the claim that Calvin initiated a dramatic shift in how knowledge is conceived, which enabled the mathematisation of uncertainty.

The third section contains the case that Calvinism was the motivation for the mathematisation of chance. It starts by arguing that probabilistic thinking before 1600 was not mathematical because it assumed chance was particular and should be addressed subjectively. It then describes the innovation of Pascal and Fermat and makes some connections to Pascal's Jansenist (Calvinist) beliefs. This is followed by a justification that Huygens' use of *spes* for mathematical expectation originated in his Calvinism. There is a discussion of the role of Puritanism, which is distinctive from Calvinism in emphasising Stoic Providentialism, in the emergence of statistics and a comment on the genesis of actuarial science. The section finishes with the counterexample of d'Alembert, who opposed Bernoulli's mathematical model in support of variolation (inoculation) based on classical *phronesis*. In the Conclusion, remarks are made on the relevance of the variolation debate to contemporary issues with vaccinations as well as concerns relating to the application of machine learning.

## A review of the evolution of Faith, Hope and Prudence

### Pre-Christian Concepts

Modern understanding is that hope is a positive emotion. The pre-Socratic Greeks were more ambivalent. Elpis (hope) was a daemon who brought feelings of comfort, but these could be misguided. This is first observed by Hesiod in the Pandora myth that describes how Zeus exacted revenge on Prometheus ("forethought") by creating Pandora ("all gifts") who opened a box that released all the ills on the world leaving only *elpis*/hope. Theognis of Megara presented Elpis as the only deity that never left the human world for Olympus (Theognis, 1931, Book 1, 1135-1150). Euripides said "Hope tricks people. Do not trust it!" (Euripides, 2010, l. 480) while in Aeschylus' *Prometheus Bound* Prometheus explained his fate with "I caused mortals to cease foreseeing their doom. … I caused blind hopes to dwell within their breasts" (Aeschylus, 1926, pp. 249-252).

Elpis was not an important Greek deity, though she was often associated with the more significant Nemesis. Nemesis' principal role was in correcting the random actions of the deity Tyche (Kershaw, 1986, 4/1-8). Tyche was a central motivating power in pre-Socratic Greece and was neither inherently good nor bad but was responsible for events over which humans had no control. It was Nemesis' role to correct any injustices resulting from Tyche's actions, and Elpis connected the individual to Nemesis. While Tyche was the motivating power, there was no point in worshipping her because she was untrustworthy (Morgan, 2015, p. 137).

The Greek daemon associated with faith, Pisitis does not appear in Homer, though *pistos* is common and was related to an aristocratic ethos (Di Somma, 2018, p. 32) (Faraguna, 2012, p. 356). *The Iliad* focuses on public *pistos* among warriors while domestic *pistos* is central to *The Odyssey*. There was a shift from the 'aristocratic' *pistos* of Homer to the 'democratic' *pistis* of the monetised polis in the mid-sixth century BCE. Money was explained as a response to the fact that individuals cannot survive on their own. If the community was in good order, the daemon Pisitis would enter it, enabling money and commerce, resolving the problem of individual limitations. If there was no trust (*apistia*), Pistis left, money was hoarded, commerce ceased, and disorder arose. By the fourth century BCE, *pistis* was being used in the sense of 'credit' and 'credit-worthiness' and was seen as the most important characteristic a merchant (Faraguna, 2012, pp. 361-362). However, wealth was not a good indicator of *pistis* since it was founded on *tyche* and not in in excellence (*arete*), justice (*dike*) or wisdom (*sophia*) (Faraguna, 2012, p. 358).

Latin culture took on many of the Greek aspects of *elpis* and *pistis*. The Romans continued to regard *spes* as ambivalent ('*spes inanis*') but they were generally more positive, with Cicero defining "*spes est expectatio rerum bonarum*" ("hope is the expectation of benefit") (Walsh, 1974, p. 34) (Clark, 1983, pp. 85-87). Spes was associated with Fortuna and hopes for prosperity with parents worshiping Spes in relation to the hopes they had for children. While *spes*/*elpis* became connected to a personal response to *tyche*/*fortuna*, *pistis*, as *fides*, related to communal mitigation of *tyche*/*fortuna* and became a fundamental feature of Roman citizenship (Morgan, 2015, pp. 107-108). The need to employ material contracts highlighted the fragility of *pistis*/*fides*, which, like *elpis*/*spes*, could be ambivalent. A recurring aspect of *pistis*/*fides* was that it related to experience and empirical evidence, which could be counterfeited and so could be misleading. (Momigliano, 1984) (Hay, 1989) (Faraguna, 2012) (Morgan, 2015).

Whereas the pagans saw *tyche* as a motivating power that was too capricious to rely on, the Jews placed hope/*qavah* in a covenant with Yahweh, their motivating power. This hope was often strained by the fact that the Jews suffered, despite adhering to the covenant. This resulted in the problem of theodicy, justifying God's omnipotence in the presence of evil and suffering, and required faith/*emuwnah*, in the covenant.

While Yahweh was the foundation of wisdom for the Jews, the pagans developed philosophy as a rational basis for managing the misfortune created by Tyche (Nussbaum, 2001, pp. xvii-xviii). A core component of this approach was the idea of *phronesis*, central to Aristotle's moral philosophy. *Phronesis* was insight born out of experience that related to particular situations, rather than universals, which were understood through *episteme* (Aristotle, 2011, VI), (Aristotle, 1926, 1366b13), (Nussbaum, 2001, pp. 294-297), (Long, 2006, p. 162), (Randall, 2011, p. 209). Aristotle highlighted the difference between *phronesis* and *episteme* with the example of mathematics. Because mathematics dealt exclusively with abstract universals, known through *episteme*, it could be understood by the young through instruction. But the young could not possess *phronesis* as they lacked experience of particular circumstances (Aristotle, 2011, VI.8).

Pagan life was influenced more by Stoic and Epicurean ideas than Aristotle. The Stoics believed that order in the universe came about through divinely constituted laws of nature. The role of reason was to come to understand those laws in order to flourish. This led to the belief that there were moral principles that were rooted in nature, *physis*, which transcended custom or opinion, and it was through reason that the natural laws (*nomos physeos*) could be understood. Cicero observed that these natural laws meant that nature was predictable and resulted in the ideas of providence and *prudentia* (Prudence). Cicero claimed that *prudentia* was built on *memoria* (memory), *intelligentia* (intelligence) and *providentia* (foresight). He focussed on *prudentia* as the virtue necessary for good government, rather than as a personal virtue, such as *phronesis*. (Annas, 1995, pp. 247-249)

The Epicureans rejected teleological explanations and regarded the natural world as being aimless (*matēn*) and this word became associated, along with *tyche*, with chance (Long, 2006, p. 160) (Reardon, 1975, p. 519). They accounted for the apparent order of the universe by asserting that while individual atoms moved randomly, in aggregate they acted in a coherent manner; a modern corollary would be statistical mechanics. They became associated with atheism because they observed evil and argued that God "either wishes to take away evils, and is unable; or He is able, and is unwilling; or He is neither willing nor able" (Lactantius, n.d., Ch XIII).

*Catholic Concepts*

Christianity emerged as a synthesis and development of Pagan and Jewish ideas. It replaced the personal *eudaimonia* of the pagans with a universal aim of salvation, 'everlasting life', granted by a transcendent God. On this basis, Christian Hope was different to Jewish *qavah* or pagan *elpis*/*spes*, which both related to mundane affairs. This shift to the spiritual resulted in the Christians adopting an unambiguously positive understanding of *spes*/hope, with Peter the Lombard (c 1096-1160) identifying Hope/*Spes* as *certa exspectatio futurae beatitudinis* ('certain expectation of future happiness') that came in salvation.

In committing a sin, Christians believe that an individual acquires the liabilities of guilt, which prevent salvation and so need to be redeemed through God, a process known as 'justification'. In the fourth century, Pelagius preached that an individual could, through their own free will, choose to be sinless. However, this was impractical for most people and the Pelagian approach undermined the role of God, and Christ, in delivering salvation. Augustine refuted Pelagius by arguing that God bestowed Grace on an individual giving them the ability to avoid sin and keeping them on a path to salvation. Developing Cicero, Augustine argued Christians avoided sin by employing *prudentia*, which required knowledge of fundamental moral principles (*intelligentia*), understanding of the world as it is (*memoria*) and an ability to plan actions towards the attainment of good (*providentia*). This *prudentia* was motivated by Grace. Augustine incorporated Stoicism by using *providentia* as the foundation for Christian Divine Providence that described God's ordering of the universe that enabled foresight founded on knowledge (Verbeke, 1983) (Ebbesen, 2004) (Ingham, 2016).

The Catholic Church came to teach that justification came from both faith and works (James 2:24) and established the doctrine that the guilt of sin could be removed through the Church, penitential acts or through purification in purgatory, an intermediate state between physical death and heaven that had emerged at the start of the second millennium (Le Goff, 1990). A consequence of the synthesis of the doctrines of penance and purgatory was the emergence of Indulgences. Indulgences were sold to Christians and, drawing on the infinite 'Treasury of Merits' made available to the Catholic Church through Christ's crucifixion, would clear a debt of sin, removing the necessity of penance or time in purgatory.

Pagan treatments of *fides*/*pistis* had centred on developing a, subjective, confidence while *elpis*/*spes* was a response to the vagaries of *tyche*/*fortuna*. Thomas Aquinas set out the Scholastic position and argued that Christian Faith and Hope focused on a transcendent, objective, and consistent God (Aquinas, 2017, II-II.17.6). Aquinas argued that Faith and Hope were appetites towards the good, rather than an intellectual process, and drew a clear distinction between Christian and more cognitive, pagan, approaches to faith. With reference to Peter the Lombard's definition "Hope is the certain expectation of future happiness" Aquinas argued that, as a God-given virtue, Hope is indubitable. However, an individual's salvation was contingent on their free-will and so there is a, subjective, uncertainty in the certainty of Hope (Aquinas, 2017, II-II.18.4).

In comparison to Faith and Hope, Scholastic *prudentia* was a *techne*, an individual's subjective application of universal, theological, principles to particular situations, and so was more Stoic than Aristotelian (Aquinas, 2017, II-II.Q47) (den Uyl, 1991, p. 103). Mundane Prudence maintained its position as the primary cardinal virtue but, along with all the cardinal virtues, it was subordinate to the Christian virtues of Charity, which unites the Christian to God; Faith, the source of theological truth; and Hope, the trust that God will deliver the Sovereign Good (Aquinas, 2017, II-II.17.6). Scholastic philosophy maintained Divine Providence as God's *prudentia* (Aquinas, 2017, I.Q22) that created and sustained the physical world; established natural law as the basis of morality; and answered individuals' prayers (Aquinas, 2017, III.64).

Aquinas also addressed the Epicurean argument that the existence of evil was incompatible with the Christian God (Aquinas, 2017, I.Q2.3). Following Aristotle, he argued that evil was the absence of God and resulted from a lack of Faith guiding reason (Aquinas, 2017, I.Q49). Following Augustine, he observed that God's concern was Divine Providence and if nature delivered misfortune, it was up to individuals to employ their Faith and *prudentia* to make the best of circumstances (Aquinas, 2005, III.64). Hence, Divine Providence was not inconsistent with evil nor with contingency in nature. In addition, Divine Providence did not inhibit free will (Aquinas, 2005, III.71-74), though God was omniscient and, if necessary, Divine Providence would always prevail through miracles (Aquinas, 2005, III.75-77).

### *Calvin's Reformation*

The aggressive sale of Indulgences, particularly in Germany, to fund the Papacy was an immediate cause of the Reformation. To counter the practices around Indulgences, the Reformation repudiated the doctrine of salvation through works in favour of Paul's emphasis on Faith (Ephesians 2:8-9), which had been developed in Augustine's rebuttal of Pelagian salvation through free will. Luther regarded the doctrine of justification by Faith alone, *sola fide*, as "The article with and by which the church stands". Calvin (1509-1564) adhered to a stricter doctrine than Luther, that of 'double predestination' where both salvation and damnation were predestined and there was no way of changing that certain fact. Luther did not develop the Catholic doctrine of Divine Providence but Zwingli, a Swiss Reformer, did and adopted a conventional Stoic perspective. Calvin, however, did create a new understanding of Divine Providence because he wanted to escape the determinism of the Stoics and Zwingli (Barth, 1960, III3 14) (Ong, 2003, pp. 15-23).

Calvin needed to innovate regarding Providence because 'double predestination' appeared to remove free-will. This was an issue of concern to Calvin before his break with Catholicism, in 1532 Calvin had observed

> The Stoics, who attribute the superintendence of human affairs to the gods, assert providence, and leave nothing to mere chance. The Epicureans, although they do not deny the existence of the gods, do the closest thing to it: they imagine the gods to be pleasure-loving, idle, not caring for mortals, lest anything detract from their pleasures; they deride Stoic providence as a prophesying old woman. (Calvin, 1969, Book I Chapter 1.2)

The Stoic position, of a fully determined universe, meant humans had no effective role to play in it while the Epicurean position regarded the gods as disinterested in human affairs and so there was no meaning to life (Reardon, 1975, pp. 528-529). Calvin sought to develop a Christian theology that explained the uncertainty of life and was distinctive from the classical, humanist, positions (Reardon, 1975, pp. 521-522) (de Petris, 2008, pp. 44-48). Calvin did this by separating his doctrines of predestination and Divine Providence.

Predestination related to the particular circumstances of an individual, it was subjective (Muller, 1986, pp. 19-23) (Ong, 2003, pp. 18-23). The Elect, a small minority of the population destined for salvation (Matt 7:14, 22:14) (Calvin, 1846, III.24.viii), were bestowed with Grace and so infused with the Christian virtues of Faith, Hope and Charity. Possessing the virtues meant the Elect were disposed to good but would, none the less, sin on account of people's innate nature after the Fall. Reprobates, the vast majority, were 'passed over' by God (preterition) and barred salvation (Calvin, 1846, III.23.i). Lacking Grace and the Christian virtues, the Reprobates would only make a good decision by accident. An individual's salvation or damnation was fixed and certain, but what happened to the individual in life depended on the decisions of themselves and others making it impossible to perceive the ultimate destination (Synod of Dort, 1619) (Reardon, 1975, p. 530) (Steinmetz, 2010, pp. 235-246). Salvation was certain but indiscernible.

If a person's behaviour has no effect on their salvation, then it is reasonable for that person to become hedonistic and focused on satisfying mundane desires. However, Epicurus advocated katastematic pleasure, based on an absence of disrupting desires (*ataraxia*), as the ultimate good. In the context of Calvinism, the Reprobates, lacking Grace, are driven by physical desires that need satisfying, the Elect, bestowed with Grace, exist in katastematic pleasure and are naturally drawn towards the good.

In contrast to the subjective nature of individual predestination, Calvin's Divine Providence was concerned with God's direction of nature and history on an aggregate level, it was universal. Calvin refused to constrain God and distinguished Divine Providence from Stoic teachings on natural law and neo-Platonic attitudes to God that regarded the universe as a God-created machine (Calvin, 1846, I.16.ii-vi) (Brunner, 1949, pp. 296,321 ff) (Barth, 1960, III3 14) (Reardon, 1975).

The separation of individual predestination and universal Divine Providence highlighted the question of theodicy, since the Elect would inevitably suffer as a consequence of their own misjudgements and the actions of Reprobates. Because God should be omnipotent, Calvin's theology attributed evil directly to God, rather than it being explained by the absence of God. Calvin justified God's creation of evil in two ways (de Petris, 2008, pp. 150-187). Firstly, Calvin noted that an act had no moral quality since morality is dependent on the intentions of the agent. This was close to the position Machiavelli had taken a generation earlier (Machiavelli, 1909, XVIII), (Randall, 2011, pp. 209-210). Furthermore, Calvin asserted that God was good, with an incomprehensible love for the Elect, but hatred for the Reprobates. Therefore, if the Elect experienced evil it was because God was using the suffering for a just purpose, such as a retributive punishment or an educative test of faith (de Petris, 2008, pp. 189-231).

This explained evil but did not securely resolve the problem of theodicy and so Calvin employed a second line based on the doctrine that God was incomprehensible, *Deus Absconditus*, "The hiddenness of God" (Isaiah 45:15). Calvin distinguished between God's revealed justice, which might be manifest in humanity's persistent experience of evil and suffering, and God's hidden justice, which was related to Divine Providence. This hidden justice was not apparent in mundane experience but relevant to an infinite God who existed in a separate dimension.

On this basis, Calvin argued that what people see as chance was a combination of Divine Providence and God's hidden justice

> that which is vulgarly called Fortune, is also regulated by a hidden order, and what we call Chance is nothing else than that the reason and cause of which is secret.
> (Calvin, 1846, I.16.viii)

Calvin accepted the random events could be coherent with Divine Providence, just as the Epicureans recognised that the random movements of atoms could be meaningful in aggregate (Partree, 2005, pp. 97-104). This innovative approach was probably rooted in neo-Epicureanism that Calvin had been exposed to in his legal training (Pitkin, 2016) (Strohm, 2019, pp. 48-51) and the Epicurean idea that the gods had interests other than the mundane.

In effect, (apparent) chance could not be taken out of human life. This was a radical shift from Stoic and Scholastic approaches to chance (Brunner, 1949, pp. 296,321 ff) (Torrance Kirby, 2003). Whereas the Scholastic approach to Divine Providence enabled foresight founded on knowledge, Calvin denied this, because Providence was indiscernible. On this basis, Calvin was consistent in dismissing both rationalism and mysticism as vain attempts to understand God's hidden justice and he called for Christians to trust in God and not to attempt to control *Fortuna.* (Reardon, 1975, pp. 527-528) (Bellhouse, 1988, pp. 68-69) (Torrance Kirby, 2003) (Gootje, 2005)

The doctrine of *Deus Absconditus*, that Divine Providence was incomprehensible, was founded on Calvin's innovations with respect to the fundamental nature of knowledge. Calvin began the *Institutes of the Christian Religion* with the assertion that "Our wisdom, in so far as it ought to be deemed true and solid Wisdom, consists almost entirely of two parts: the knowledge of God and of ourselves" (Calvin, 1846, I.1.i), with knowledge being founded on either faith or experience (Partree, 2005, pp. 29-42). Developing Paul's scripture (e.g. Col. 1:9), Calvin rejected the idea of diverse types of knowledge – *sapientia*, *scientia*, *intellegitia*, *prudentia* – synthesising them all into a single, essential, knowledge of God (Essary, 2017, p. 109). This synthesis is a significant departure from Aristotle, who emphasised that *phronesis*, originating in personal experience, was distinct from the universality of *episteme* or *sophia*. The distinction was present in the Stoic's understanding of the difference between *prudentia* and *scientia* and the Scholastics distinction of appetitive and intellectual virtues that regarded *prudentia* as a subjective response to events. However, this distinction was superfluous if the vagaries of life were "nothing else" than God's hidden reason (Calvin, 1846, Aphorism 13).

Calvin defined the relationship between Faith and Hope for the Elect in the context of his doctrine as

> Thus, faith believes that God is true; hope expects[2] that in due season he will manifest his truth. Faith believes that he is our Father; hope expects that he will always act the part of a Father towards us. Faith believes that eternal life has been given to us; hope expects that it will one day be revealed. Faith is the foundation on which hope rests; hope nourishes and sustains faith. (Calvin, 1846, III.2.xlii)

Calvin adopted the position that Hope, bestowed on the Elect, relates to the Elect's certain salvation, provides confidence and is a kind of waiting. In this, well-known, passage, Calvin highlighted the cognitive, rather than appetitive, nature of Hope and denied the characteristics that Aquinas had attributed to Hope (Aquinas, 2017, II-I.40.2). Whereas Aquinas saw an individual's free will inhibiting their salvation (Aquinas, 2017, pp. II-II.18.4) this is not possible in Calvin's theology, founded on predestination, and so Christian Hope is indubitably "the certain expectation of future happiness" (Synod of Dort, 1619, I.13) (Bromiley, 1953). While Hope is certain for the Elect it is non-existent for the Reprobate, such as a contemporary atheist, who is not bestowed with Grace, has no Faith, does not Hope for salvation and so has no anxiety about salvation. Similarly, Calvin argued against Scholastic understanding of 'blind' Faith with "Faith consists not in ignorance, but in knowledge" (Calvin, 1846, III.2.ii) and "that it is a firm and sure knowledge" (Calvin, 1846, III.2.vii). Here, 'knowledge' is knowledge of God, which is, necessarily, certain for the Elect.

This doctrine revolutionised the nature of Faith and Hope. Pagan *pistis*/*fides* came out of personal experience while *elpis*/*spes* related to individual circumstances. Catholic and Lutheran Faith and Hope, while directed to the universal and objective aim of salvation, were similarly subjective because the individual played a part in their own salvation. In these circumstances, Hope was doubtful. Calvinism rejected the idea that there could be an "uncertain certainty" (Synod of Dort, 1619, I.R7) and considered Hope and Faith, intellectual virtues, as being universally experienced by the Elect and the accompanying salvation was certain and immutable (Synod of Dort, 1619). Uncertainty in the certainty of God and salvation was subjective, in the minds of people, rather than a truth ordained by God. Because the relationship of individual experience to certainty is changed, the distinction between individual *phronesis* and universal *episteme* and *sophia* is lost in Calvin's theology.

---

[2] "hope expects" is "*espérance attend*" in French (Calvin, 1859) and "*spes expectat*" in the 1536 Latin edition (Calvin, 1536).

*Experimental Calvinism and Puritanism*

The Calvinist doctrine of predestination is harsh, not least because only a small minority would be Elected[3]. Despite the objective certainty of Hope and Faith, the reality was that, because of God's hidden nature, Calvinists could not be confident that they were truly experiencing the katastematic pleasure that came with Grace and delivered certain Faith and Hope. This resulted in anxiety (Synod of Dort, 1619, I.16) (de Petris, 2008, pp. 232-279) that was a key feature of 'Experimental Calvinism', often identified with Puritanism, which developed Calvin's attention to experimental knowledge and had emerged in Britain in the early seventeenth century and spread across Europe and then to North America.

While the Dutch Reformed Church was constituted on doctrine laid out in the *Canons of Dort* (1619), Puritan doctrine was presented in the *Westminster Confession of Faith and Catechisms* (Assembly of Divines, 2005) developed by English and Scottish Protestants between 1643 and 1647. The Dutch Canon emphasised the significance of Election, with Faith following, and that God is hidden and incomprehensible. The Westminster Confession focused on the 'Covenant of Grace' that enabled Faith in the Covenant that was supported by the evidence of 'Special Providence', relating to God's direct intervention in the world and individuals' lives. The significance of Providence in the Westminster Confession, while it is absent from the Canons of Dort, is that the Epicurean aspects of Calvin's 'hidden' God are being diminished with the Puritan God being more involved in the lives of the faithful, suggesting the influence of Stoicism from Scholastic and Humanist ideas.

The indubitability of Election/Reprobation implied to the Puritans that salvation could be considered objectively, in contrast to Catholic salvation, which was dependent on individual actions and so was subjective. One feature of Experimental Calvinism was the scrutiny of daily events to discern evidence of God's involvement in the lives of the Puritan that signified Election (Oxenham, 2000), (Cambers, 2007, p. 798), (Reeves, 2015, p. 613), (Ryrie, 2016, p. 51). This appears to be contradictory to Calvin's insistence that God's will was hidden and, again, signifies the resurgence of Stoic order into Reformed belief. One habit of Experimental Calvinism was the practice of keeping detailed 'spiritual diaries', or 'spiritual autobiographies', that recorded daily events and emotions. These are now widely regarded as being instrumental in the development of the English novel, in particular with Bunyan (1628-1688) (Stachniewski, 2008) and Defoe (1660-1731) (Starr, 1965), and their influence on the development of modern, experimental, science has recently been discussed by McKeveley (2018). Given the impact of these practices outside mathematics, it is reasonable to believe they would have also affected mathematics.

---

[3] Calvin (1846, p. III.25.xii) appears to argue that the idea of Hell is figurative and on death the Reprobate is simply extinguished ('conditional immortality') and so is not condemned to everlasting suffering. The meaning of Calvin's comments is contested.

*The Evolution of Phronesis and Prudence After Calvin*

Machiavelli, some forty years older than Calvin, regarded *phronesis* as the primary virtue, in both *The Prince* and the more substantial *Discourses on Livy*. *The Prince* (1513/1532) discussed how the ruler of a newly acquired state could secure stability (Randall, 2011, p. 209), the *teleos* of a state. Machiavelli argued that a prince should employ *phronesis* in support of this *teleos* and that actions should be judged in relation to this *teleos*, not morality (Machiavelli, 1909, Chapter XVIII). This was coherent with Aristotle's *phronesis* (Randall, 2011, pp. 209-210) but at odds with Catholic doctrine *prudentia* being the application of universal moral principles. Machiavelli used history to deliver the experience (*memoria*), which forms the basis of political reasoning. He emphasised that no two situations were the same and so there can be no certain guidance on action. In particular circumstances, *virtu*, or ability (*intellegentia*), needed to be employed to determine the specific actions in a deliberative manner (Cox, 1997).

Montaigne's *Essais* (1570-1592) was concerned with making moral judgement in the face of Fortuna where opinions, not substantive truth, dominated. This suggests an identification with *phronesis* but, in a comment on a story from Seneca's *De clementia*, Montaigne dismissed prudence with "vain and futile a thing is human prudence; throughout all our projects, counsels and precautions, Fortune will still be mistress of events" (Montaigne, 1877, XXIII). In 1581, the Catholic censor recommended that Montaigne remove the frequent references to Fortuna and suggested replacing them with a more theologically sound Providence. Montaigne did not replace Fortuna with Providence, indicating that he considered the two distinct, but in subsequent editions of *Essais*, Providence became more prominent. Montaigne, writing after the emergence of Reformed religion in France, took the position that Divine Providence was general and inaccessible to humans, who were left experiencing capricious fortune. The significance of these observations is that while Montaigne shared some of Calvin's attitudes, he appears to have distinguished chance from Providence. (Legros, 2009) (Regosin, 2009)

Descartes (1596-1650) was inspired to write the *Discourse on the Method* (1637) as a reaction against Jesuit teaching on probabilism, the attitude that the most convenient opinion was justified providing that it was supported by some authority (MacDonald, 2002). Descartes was dissatisfied with this ambiguity and left France for the Calvinist dominated Netherlands, where mathematics, what the Dutch called *wiskunde* 'the art of certain knowledge', was ascendant. In the *Discourse*, Descartes focused on knowledge of truth, *scientia*, gained through natural reason and argued never to accept anything as true unless you were certain that it was true (Descartes, 2008, Part II). To escape radical scepticism, he adopted an Augustinian argument (Augustine, 2000, Book XI, 26) that the fact that he doubted meant that he existed, and this certainty was the root of all knowledge (Descartes, 2008, Part IV). From this notion, Descartes constructed a theory of knowledge using the deductive method of mathematics to ensure its veracity (Descartes, 2008, Part I).

Descartes struggled to apply his philosophy coherently to the problem of chance. In correspondence to (the Calvinist) Elisabeth of the Palatinate in 1645, he took a Stoic position that the role of reason was to come to terms with misfortune. Later, in the 1649 *Les Passions de l'âme* (The Passions of the Soul), Descartes initially argued that reason and forethought could manage misfortune, but at the end of the book he concludes that they were of little help when dealing with the particular realities of an uncertain world (Gilby, 2009). Descartes was interested in addressing the problem of uncertainty but was unable to develop a conceptual framework to do so.

Hobbes' (1588-1679) principal motivation was to mitigate the turmoil caused by the Puritan rebellion against the catholic king in England (Shulman, 1988). He rejected the classical framework of morality, as a rational process aimed at individual flourishing, in favour of a framework that was based on individual desire to live in a society, a Commonwealth, that would protect them (Hobbes, 2017, Chapter XVII) (den Uyl, 1991, p. 120). Hobbes sought an objective, immutable, morality that could be defined in terms of the Commonwealth (den Uyl, 1991, p. 114) and he rejected the idea that there could be a personal morality. An accomplished mathematician, Hobbes employed mathematical arguments in favour of his political science (Hobbes, 2017, Chapter V) (Randall, 2011, pp. 212-213) and Hobbes sought a scientific (*episteme*) morality founded on universals, relevant to the Commonwealth, rather than a reasoned morality (*phronesis*) reacting to particular circumstances.

In *The Elements of Law* (1640) and *De cive* (1642), Hobbes explicitly rejected prudence as merely "opinion"; it would only ever be correct by chance. Prudence was a virtue only in so much as a synthesis of qualities such as good memory, focus and intelligence that resulted in sound judgement (Hobbes, 2017, Chapters VIII, XII). In addressing particulars, prudence was probabilistic, not in the mathematical sense but founded on opinion, and so could not be related to certain knowledge (Hobbes, 2017, Chapter XLVI). Hobbes, like Descartes, is a mathematician seeking certainty in the face of apparent uncertainty and, on this basis, dismisses prudence, which had been the paramount virtue in both Classical and Scholastic philosophy.

Spinoza (1632-1677) was an apostate Jew from Amsterdam, a part of the Reformed Netherlands that was not strictly Calvinist. His most influential work, *Ethics*, was published posthumously in 1677. Spinoza echoed Plato, Augustine and Descartes in arguing that mathematics provided the means of discerning truth (Spinoza, 1901, I Appendix 2-3) and the argument of *Ethics* followed a mathematical format presenting a deductive chain that proved propositions out of definitions and axioms. Spinoza's most famous quote "*Deus sive natura*", "God or nature", summarised that there was no distinction between mind, matter and God. People believed themselves to have autonomy because, being finite, they could not comprehend the infinity of God. On this basis, Spinoza argued that the purpose of the individual was to lift themselves out of a mundane perspective so that they could understand the totality of creation, coming to understand the true nature of God's will, which were manifest in the laws of nature. In Spinoza there is a unification of *episteme*, *sophia* and *phronesis* but the possibility of understanding God with the absence of Calvin's mundane experience of chance, highlighting Spinoza's Stoicism (DeBrabander, 2007) (Miller, 2015, pp. 1-6). This is exemplified in Spinoza's belief that events are necessarily determined "Nothing in the universe is contingent, but all things are conditioned to exist and operate in a particular manner by the necessity of the divine nature." (Spinoza, 1901, I.P29).

## Calvin's theology and the development of the mathematics of chance

*The Scholastic Contribution – Olivi*

Kaye (1998, pp. 120-122) and Franklin (2001, pp. 265-269) make the case that the origin of mathematical probability is in the spiritual Franciscan Pierre-Jean Olivi's analysis of Justice in commerce. These arguments are in the context of the thesis that the mathematisation of science emerged out of the medieval monetisation of society (Restivo, 1982, p. 128), (Hadden, 1994, p. 84), (Crosby, 1997, pp. 69-74), (Goetzman, 2016, pp. 203-275), (Kaye, 1998).

The Kaye/Franklin argument originates in Aquinas examining the morality of commercial exchange where, using an example from Stoic philosophy, Aquinas (2017, II-II.77.3) argued that uncertainty in events could allow for profits. Olivi disagreed with Aquinas and argued that the Dominican's analysis demonstrated *prudentia* but not Christian Charity and Faith which gave individuals the ability to make a, subjective, judgement defining a 'just price' (Zimmermann, 1996, pp. 264-268) (Kaye, 1998, p. 25) (Monsalve, 2014).

According to Kaye and Franklin, the Scholastics' argument implies a significant conceptual leap in that chance, an immaterial concept, was being quantified since the expectations were being expressed as prices (Kaye, 1998, pp. 119-124) (Franklin, 2001, pp. 265-267). However, the quantification of beliefs through their expression as prices was not novel. Indulgences satisfied a *pena* for sin, just as older Roman Law allowed all crimes, apart from murder, to be satisfied by fines. The concrete, quantified payments relating to abstract concepts, such as sin and injury, was not novel. Moreover, the Scholastic understanding of probable opinion was that it was subjective and motivated by individual Grace (Kaye, 1998, pp. 154-155,168-169). There is no evidence that Olivi believed that the prices could be established through *scientia*, rather than decided based on subjective *prudentia*, guided by individual Grace, so he was not concerned with an objective mathematical idea of probability. When a modern, small scale, business estimates its future income and expenditure it is done based on personal experience informing judgement, classical *phronesis*/*prudentia*, rather than on the basis on the mathematical idea of expectation. The existence of contracts that relate to uncertain values does not imply mathematical conceptions of chance, probability or expectation.

Scholasticism was based on both theology and jurisprudence While the theologians regarded Faith and Hope as appetitive powers related to God, the Roman idea of *spes* as a cognitive power related to *fortuna* was preserved in the medieval legal treatment of aleatory contracts. In legal treatments, particularly in relation to inheritances, there were references to a *spei,* the genitive form of *spes,* and usually translated as an 'expectancy' (Zimmermann, 1996, pp. 245-249) (Franklin, 2001, pp. 258-261). Franklin (2001, pp. 348-352) suggests that Huygens' use of *spes* in relation to probability was motivated by this term *spei* and argues that mathematical probability evolved out of Scholastic jurisprudence. The idea of *spei* exists in current English and American law, in the tort of an interference with an expectancy/expectation of inheritance. These expectations are quantified but this is not done with reference to mathematics. Dickens' novel *Great Expectations* does not suggest mathematics, nor does the discussion of 'expectations' of the results of court cases in *Bleak House.* Rather, they relate to the Roman conception of *spes*, which was highly subjective. Franklin highlights the particularity of legal judgements and their uncertainty (Franklin, 2001, p. 350) and so undermines the thesis that legal thinking can be the origin of a universal, mathematical, probability.

Franklin supports the legal basis of probability by highlighting the pioneers of probability who had a training in law. However, Fermat was the only one whom Franklin identifies who actually practised law. The significance of legal training is in that law students, notably Calvin, were exposed to Epicureanism is recognised in the argument presented here (Pitkin, 2016) (Strohm, 2019, pp. 48-51).

*The Humanist Contribution – Cardano*

Cardano's mid-sixteenth century work on probability, *Liber de Ludo Aleae*, considered the quantification of chance by investigating the *episteme* of gambling. Cardano was a humanist interested in understanding the behaviour of people, with different needs and making different choices, and, given that variety, to find a viable way towards individual happiness (Canziani, 1992, pp. 295-296). He mirrored Machiavelli's innovation that *phronesis* did not always lead to Christian morality with the idea that *episteme*, similarly, did not necessarily result in moral goodness, breaking with both Catholic and classical philosophy that assumed knowledge was intrinsically good. This idea was synthesised with a conventional Stoic line that through *sophia/sapienta* people could navigate the ills of the world (Cape, 2003, p. 40) (Bracali, 2008).

In the *Liber,* Cardano observed that the fundamental principle of gambling was equal conditions, which was associated with Aristotelian concepts of Justice (Bellhouse, 2005), and argued that, on account of symmetry, the chances of throwing one, three or five are the same as throwing a two, four or six. Cardano finished this point with the well-known observation "these facts contribute a great deal to understanding but hardly anything to practical play" (Cardano, 1953, p. 194). This observation exemplified the classical distinction between *episteme*, knowledge of universal principals that affect the chances, which Cardano addressed mathematically, and *phronesis*, practical wisdom necessary for sound judgement. Cardano still accepted that judgement was founded on *phronesis/prudentia*, his innovation was in presenting an approach to *tyche/fortuna* that is *episteme/scientia* but he could not make the conceptual leap that science could be employed to manage chance.

*Indirect Calvinism – Pascal and Fermat*

The canonical birth of mathematical expectation, a universal means of calculating values used for decision making, is in the Pascal-Fermat correspondence that solved the *Problem of Points*. The *Problem* can be summarised as follows

> Two players, *F* and *P*, are playing a game based on a sequence of rounds, and each round consists of, for example, the tossing of a fair coin. The winner of the game is the player who is the first to win 7 rounds, and they will win 80 francs.

The *Problem of Points* is how the 80 francs should be split if the game is forced to end after *P* had won 5 rounds while *F* had won 4. It originated in the *abaco* tradition of using 'stories' to give examples of how to solve problems in commercial calculations (Sylla, 2006).

Pacioli had offered a solution that was statistical, that the pot should be split 5:4 reflecting the history, *memoria,* of the game. Cardano had argued that this was absurd since it would give an unfair result if the game ended after one round out of a hundred or when *F* had 99 wins out of a hundred to *P*'s 90. Cardano made the point that a better solution would be arrived at by considering what would happen in the future, *providentia*, not the past. To do this, one had to account for what 'paths' the game would follow. Despite this insight, Cardano's solution is generally regarded as insufficient, and an improved solution was provided by Pascal and Fermat in a series of letters to each other in 1654.

Pascal or Fermat, it is not certain who, realised that when Cardano calculated that *P* could win if the game followed the path *PP* (i.e. *P* wins and *P* wins again), this represented four paths, *PPPP*, *PPPF*, *PPFP*, *PPFF*, of the game. It was the players' 'choice' that the game ended after *PP*, not a feature of the game itself. This is a mathematical insight, abstracting from the particularities of the game to identify the fundamental structure of the problem and justified the view that Pascal and Fermat were the first to develop a mathematical theory of probability. Having discerned the underlying structure, calculating the proportion of winning paths came down to applying the Arithmetic, or Pascal's, Triangle resulting in the Binomial Model.

The problem is binary, just as Election/Reprobation is binary, and the players are seeking to know their terminal state based on immediate experience, which suggests Experimental Calvinism. Furthermore, by considering parts of the paths that are not explicit but are fundamental to the overall structure the solution relates to the Calvinist idea of God's hidden justice. Pascal's family were Jansenists and followed the doctrine of double predestination close to Calvinism. Fermat was a Catholic, however he had an affinity with Reformed religion having had his legal training at Orléans, which had been Calvin's *alma mater* and attracted Calvinist students from across Europe, and he had chosen to practice law in Castres, a Huguenot city where his only clients would have been Calvinists (Chabbert, 1967), (Barner, 2001, p. 14). It can be inferred that both Fermat and Pascal would have been familiar with the Reformed doctrine that the ultimate ends were certain though indiscernible, even if Fermat did not personally espouse it.

While corresponding with Fermat, Pascal had distanced himself from Jansenism. However, after being involved in an accident a few months after solving the *Problem*, he returned to the sect and for the rest of his life he focused his efforts on supporting the Jansenist cause, primarily in his *Lettres Provenciales*, a satirical attack on the Jesuits. Pascal's *Wager* (Pascal, 1958, p. 233) was formulated at this time (1656-1662), but only published posthumously in the *Pensées*. The *Wager* was an attempt to apply reason to the question of whether to follow doctrine without having to rely on explicit evidence (Hacking, 2006, p. 64). It is an expression of the doubt that those committed to predestination had, which needed to be resolved by Faith. The argument about the evidence for the existence of God was not novel (Franklin, 2001, pp. 249-251) but it reveals Pascal's conventional Calvinism, that a person bestowed with Grace possesses Faith and certain Hope. Pascal wrote

> there is sufficient evidence to condemn, and insufficient to convince; so that it appears in those who follow it, that it is grace, and not reason, which makes them follow it; and in those who shun it, that it is lust, not reason, which makes them shun it (Pascal, 1958, p. 563)

Franklin comments on this as follows

> The idea that the will supplies decision where certainty is not available has here got out of hand. Where some would draw the conclusion that one ought therefore to study the probability of the evidence more carefully, Pascal concludes that there is certainty for those who seek sincerely (Franklin, 2001, p. 254)

Franklin's comment reflects the Catholic approach that Faith is an appetitive power that guides reason so that salvation is ultimately determined by an individual's free will. The role of individual will in Catholic salvation makes it impossible to predict[4] with any certainty. The Calvinist view, reflected in Pascal's writing, was that Grace bestows Faith, a cognitive power, on the Elect and, because of Grace, there is an objective certainty in salvation. However, the certainty of salvation is difficult to discern from mundane experience and the Calvinist Hopes, or expects, salvation.

---

[4] Note that the Reformation had a significant impact on Catholic practice. Before the Reformation a Catholic could immediately absolve themselves of a lifetime of sin through Indulgences or the Sacrament of Penance.

The idea that there are facts hidden behind noisy data, and by analysing the data objective expectations can be formed, is fundamental to both statistics and Experimental Calvinism. The effect of Calvin's theology on the emergence of probability and statistics is clear in Pascal's thinking. It is explicit in Pascal's Wager and implicit, in thinking in terms of paths and expectations, in the Fermat-Pascal solution to the *Problem of Points*. Both highlight the experiential anxiety Reformed Christians felt regarding Election (Synod of Dort, 1619, I.16), (de Petris, 2008, pp. 232-279), but is absent in conventional Catholicism, where salvation was contingent on an individual's actions ('works'), and, being subject to individual will, was inaccessible to science.

*Direct Calvinism – Huygens*

Christiaan Huygens was born into a prominent Dutch family in 1629. The family's relationship to the House of Orange implies that Huygens grew up in an environment committed to the Calvinist doctrine laid down in the *Canons of Dort* (Stoffele, 2006, pp. 61, n 202). Huygens was intended to follow his father into state service and so studied law, which provided the humanist education required by a courtier and included mathematics instruction under van Schooten (Stoffele, 2006).

In the second half of 1655, Huygens visited Paris and was told about the *Problem of Points*, but apparently not of its solution, which he was able to deduce independently. He returned to the Netherlands and wrote *Van Rekeningh in Spelen van Geluk* (*On the Reckoning at Games of Chance*) which would appear in van Schooten's *Exercitatonium Mathematicarum*, a university textbook, as *De Ratiociniis in Ludo Aleae*, in 1657. Huygens had to translate the Dutch text of *Van Rekeningh* into Latin for van Schooten's textbook. The book opens with

> "If I expect *a* or *b*, and have an equal chance of gaining either of them, my Expectation is worth ($a$ + $b$)/2." (Huygens, 1714)

When the text came to be translated into Latin, for what is written as 'Expectation', Huygens considered using *spes* (Hacking, 2006, pp. 93-95) and, since the eighteenth century, this option is reflected in the French use of *espérance* when referring to mathematical expectation. Van Schooten's published edition used *expectatio*, giving the English term 'expectation' as in this English translation.

Huygens would have been influenced by humanism, his legal training, and his religious upbringing. Humanist 'hope', as in Cicero's definition of *spes* as "*expectatio rerum bonarum*", was regarded as positive and to use such a loaded term for mathematical expectation seems incongruous. The more ambivalent *elpis* might have been more appropriate but there is nothing to suggest this was Huygens motivation. Had Huygens been influenced by his legal training, he would have chosen to use *spei* or *sperate* rather than *spes* (Zimmermann, 1996, pp. 245-249).

An explanation for Huygens' consideration of *spes* is in his religious upbringing and in Calvin's formula "hope expects"/"*spes expectat*" (Calvin, 1536)/"*espérance attend*" (Calvin, 1859). Calvinist doctrine was that Hope, the "certain expectation of future happiness" bestowed on the Elect through God's gift of Grace, represented a universal and indubitable idea, not the subjective hope of the pagan or Catholic contingent Hope in salvation (Aquinas, 2017, II-II.18.4). Similarly, Calvin's doctrine had changed the nature of Faith. Pagan *pistis*/*fides* was founded on personal experience just as Catholic Faith was subjective directing an individual's own reason and free will towards salvation. Calvin's Hope and Faith, experienced by the Elect, were different in that they were universal, objective and related to certainty (Synod of Dort, 1619) (Bromiley, 1953).

The argument for the Calvinist origin for the choice of *spes*/*espérance* in the context of mathematical expectation is supported by the fact that the earliest recorded use of the French technical term "*l'espérance mathématique*" was in a letter written by Gabriel Cramer in 1728, while he was chair of mathematics at the *Académie de Calvin* in Geneva. Pascal and Fermat did not use the term *espérance*, they discussed shares and values (*pistoles*) (Lee, n.d.).

The universal and immutable certainty of Hope and Faith opened the door for mathematics, universal and indubitable, to address the contingent nature of life and was coherent with both Hobbes' and Descartes' desire for objectivity. However, the changes Calvin initiated in relation to Hope and Faith proved less significant than the change to Prudence. Before Calvin, western philosophy, notably Descartes less than ten years before Huygens tackled probability (Gilby, 2009), addressed judgement and decision making through *phronesis*/*prudentia* while science, and mathematics, was addressed by *episteme*/*scientia*. Calvin dismissed the idea that there was a distinction, synthesising all types of knowledge into aspects of knowledge of God (Calvin, 1846, I.1.i). Then, he asserted that chance "is nothing else" than God's hidden reason, implying that it can be understood on the basis of having developed a universal type of knowledge. A hundred years after Calvin, mathematicians influenced by Calvinist theology were able to address the particularities of chance through the universal generalisations of mathematics, something Cardano, Descartes and Hobbes had felt unable to do.

*Experimental Calvinism*

Hacking (2006, pp. 31-48) argued that the emergence of probability was a consequence of a shift in attitudes to 'internal evidence', which he identified as having come out of changes in medical practice after Paracelsus (1493-1541; Calvin 1509-1564) that focused on reading 'signs' in nature and led to the concept of 'sign-as-evidence'. However, an alternative, and much more widespread, origin could be in Calvin's doctrine that emphasised experiential knowledge and claimed certainty in the presence of a lack of explicit evidence, the problem motivating Pascal's *Wager.*

Hacking does not mention the practices of Experimental Calvinism or the 'spiritual diaries' that emerged in response to the mundane uncertainty of Election. He does, however, come close to recognising the significance of Calvinist theology in a discussion of the influence on the emergence of these 'sign-as-evidence' of the *Port Royal Logic*, written around 1662 by two leading Jansenists, Antoine Arnauld and Pierre Nicole (Hacking, 2006, pp. 73-80). The *Logic* ends by addressing judgement of the future based on past experience (Arnauld & Nicole, 1850, pp. 358-362) and Hacking attributes this as being the first place where the term probability is used to refer to a measurable quantity (Arnauld & Nicole, 1850, pp. 360-361).

Hacking does remark on a connection between Puritanism and the genesis of modern science in the work of John Wilkins (2006, pp. 80-84). Wilkins wrote on Divine Providence, highlighting the difficulty of discerning God's intention in mundane events

> Remember we are but short-sighted, and cannot discern the various, references, and dependancies, amongst the great affairs in the world, and therefore may be easily mistaken in our opinion of them. (Wilkins, 1649, p. 72).

Hacking highlights how Wilkins began to recognise a difference between 'internal' and 'external' evidence and would initiate what Hacking calls 'Royal Society theology' that argued that the well-ordered universe was evidence of God's existence.

An associate of John Wilkins who and shared his interest in the importance of the examination of experiences in the development of 'inner sense' was the Puritan John Locke (Connolly, 2014, pp. 253-255). Locke discussed probability in Chapter 15 of Book IV of an *Essay Concerning Human Understanding* (1689) where he introduced the term 'probability' in the context of proof with the observation that certain truth is elusive in practice and so defined probability as the "likelihood of truth" (Locke, 2017b, XV 3), and connected the word probable to faith and belief.

This is significant as Locke is describing probability as something objective whereas earlier approaches focused on probability as subjective opinion. This was a developed ideas presented in the *Port Royal Logic* and other minor mathematical works of the time continued in this direction, including Craig's *Thealogiae Christianae Principia Mathematica* (*Mathematical Principles of Christian Theology*) of 1699 that looked at the credibility of human testimony and the *Calculation of the Credibility of Human Testimony*, published anonymously in Philosophical Transactions of the Royal Society in 1700.

Locke and Craig were both non-conforming Protestants and their arguments, more Stoic than Epicurean, reflected the influence of the Westminster Confession rather than the Canons of Dort. The general significance of Puritanism to mathematical probability is further highlighted by the fact that the Puritans were interested in understanding the idea of chance (Bellhouse, 1988) and considered the role of providence in gambling. The casting of lots was legitimate if seeking guidance from God, but gambling with dice was 'pure contingency' and illicit (Bellhouse, 1988, pp. 67-69), (David, 1998, pp. 13-20). This interest could be theologically justified because of Calvin's assertion that "Chance is nothing else than that the reason" of God (Calvin, 1846, I.16.viii).

*Divine Providence and the Emergence of Statistics*

The genesis of statistics was in Graunt's 1662 investigation of mortality data that revealed the stability and near equality of the sex-ratio at birth. Graunt, who had been raised a Puritan (Lewin, 2004), could not decide whether the higher male birth rate reflected the higher adult mortality of men in war or that the equality implies Christian monogamy is more natural than Islamic polygamy (Hald, 1990, p. 93) with both explanations implying Divine Providence and God's hidden justice. Lodewijk and Christiaan Huygens then engaged in correspondence on how Graunt's data could be used to accurately price life annuities, essential to the Dutch Republic's public finances. De Witt solved the problem the Huygens brothers were investigating in 1671 and Halley produced an improved lifetable based on data from Breslau/Wroclaw in 1694.

Graunt's data played a part in the most important application of mathematics in support of 'Royal Society Theology', John Arbuthnot's *An Argument for Divine Providence, Taken from the Constant Regularity Observed in the Births of Both Sexes*, presented to the Royal Society at the end of 1710. Arbuthnot was the son of a minister of the Reformed Church of Scotland and had been educated at Marischal College, which had been established to train Reformed clergy, though he was not an ardent Reformer (Ross, 2004).

The *Argument for Divine Providence* is significant because it initiated the use of mathematics to test statistical significance (Hacking, 2006, pp. 166-171), (Kemp, 2014). Arbuthnot began with Cardano's observation that symmetry implies equal chances of male and female births. He then argued that the fact that male births exceeded female births in the 82 years of Graunt's data was so unlikely as to indicate Divine Providence (Bellhouse, 1989, p. 255).

Arbuthnot was seeking to prove that chance had no role in determining the ratio of male and female births and in doing so provide a mathematical counter to Epicurean, atheistic, beliefs that chance dominated (Kemp, 2014, p. 473). Arbuthnot's argument was criticised, mathematically, by Nicolas Bernoulli in 1712 (Shoesmith, 1987) who argued in favour of Epicurean chance that if the equiprobable assumption was rejected, and the probabilities of male births set at 18/35, then chance could not be dismissed. De Moivre countered Bernoulli with the observation that the 18/35 randomiser had been made "by some (Divine) Artist" (Kemp, 2014, pp. 483-484).

These discussions concerning the quantification chance and in providing evidence for Divine Providence, as "likelihood of truth" or as an aspect of moral reasoning, are reminiscent of classical *fides*/*pistis* and Stoic natural law and suggests a genesis of mathematical probability before Calvin. However, in the seventeenth century the Calvinist idea of God's hidden justice is fundamental but lacking in earlier treatments. Moreover, quantification is not mathematisation. Before the Reformation, approaches to chance were founded in subjective *phronesis*, as has been highlighted above, while Calvin's theology had removed the distinction between subjective *phronesis/prudentia* and objective *episteme/scientia*, enabling the mathematical treatment of chance.

### *Mathematics in practical decision making – Bernoulli*

After Huygens' *De Ratiociniis*, the next important mathematical text on chance was Jacob Bernoulli's *Ars Conjectandi* (c. 1700-1705). Bernoulli was raised in the Reformed Church directed by the Helvetic Confession, which originated in the Stoic Zwingli and is closer to the Westminster Confession than Canons of Dort. Bernoulli began the final part of the *Ars* by making a novel distinction, that certainty is either objective or subjective (Bernoulli, 2005, p. 8) and made the case for probabilities being useful in moral, practical, judgements, not just for general rules in gaming. These observations are significant because subjective decision making had historically been governed by *phronesis/prudentia* and Bernoulli was categorising this as a sub-class of *episteme*/*scientia*, highlighting Calvin's unification of knowledge and the diminution of *phronesis*/*prudentia*. At the time, this approach was distinctive from Leibnitz who was arguing, in a classical phronetic manner, that practical, moral, judgements were peculiar and so not amenable to the general rules of mathematics (Meyns, 2017, pp. 230-250).

*The Reformed Church and Actuarial Science*

Nicolas Bernoulli applied his uncle Jacob's theories to insurance contracts using Graunt's life table in his thesis*, De Usu Artis Conjectandi In Jure* (*On the Use of the Art of Conjecturing in Law*), of 1709. The theory of insurance mathematics was then developed by de Moivre in the 1720s and by Thomas Simpson in the 1742. The first practical application of mathematics to manage an insurance fund, either in general insurance or in life insurance, was the Scottish Ministers' Widows' Fund, established in 1744 (Hare & Scott, 1992). The Fund was a consequence of the reformation of the Scottish Church in 1560 that allowed married men to be ministers and in 1672 the Scottish Parliament legislated that ministers' widows were entitled to stipend. By the 1730s, this one-off payment was regarded as insufficient and that the Church of Scotland should set up a 'general fund' from which pensions for widows could be drawn.

The instigator of this initiative was Robert Wallace. Wallace had been born in 1697, like Bernoulli and Arbuthnot, the son of a cleric and became the Moderator (head) of the Church of Scotland 1743-1746 (Luehrs, 1987), (Ahnert, 2015, pp. 78-93). In the 1730s Wallace had become interested in economic affairs (Peterson, 1994, p. 6) and in 1741 he directed another cleric, Webster, to gather statistics on ministers' deaths and their dependants for the period March 1722 to March 1742. Meanwhile Wallace tackled the mathematical theory behind the problem using Halley's Life Table and in 1743 he calculated the premium rates and estimated the cash-flows for the proposed fund. The Fund was established in 1744 and the fund's capital never deviated by more than 5% from estimates calculated in 1748 (Dunlop, 1992, p. 18). It was replicated in the Presbyterian Ministers Fund of Philadelphia in 1761 and the following year the Equitable Life, open to the public, was founded in England.

It is remarkable, and a fact that should not be ignored, that the first insurance fund to be managed mathematically was established by Reformed clergymen, rather than for commercial profit, noting that de Moivre had been working with merchants and bankers in London. Wallace and Webster would have been committed to the certainty that, for the Elect, Faith believes/Hope expects, and Charity is the result. On this basis, and God's hidden justice, they had objective Faith in the data, mathematical statistics, and based on certain Hope, mathematical expectation, were confident that the fund would deliver charitable ends despite the apparent randomness of individual lives.

*Universal Mathematical Models Replace Particular Wisdom – The Variolation Debate*

A clear expression of the differences between the classical approach to *phronesis/prudentia* and the Reformed approach to the use of mathematics in practical decision making comes in the debate in the 1760s between d'Alembert and Daniel Bernoulli on smallpox variolation, the practice of purposefully infecting someone with a weakened virus to prevent more serious infection.

D'Alembert was an apostate Catholic, born after the suppression of Calvinism in France (1685. He regarded mathematics relating to "intellectual abstractions" of "general and common properties" discerned out "particular and palpable properties" experienced by senses (d'Alembert, 2009), rather than the uncertain particularities of the world (Daston, 1979) (Colombo & Diamanti, 2015). On this basis, d'Alembert accepted that combinatorics were sufficient for simple, closed, games of chance but mathematics had little to say about questions of normal life, such as insurance, where the complexity of the problems meant that "only experience and observation can instruct us …, and only instruct us approximately" (Daston, 1998, p. 108 ff), (Colombo & Diamanti, 2015, p. 189). D'Alembert was making the classical case that the uncertainty of life cannot be addressed mathematically, and his position contrasts with that of his contemporary Wallace's.

On this basis, d'Alembert came into dispute with Daniel Bernoulli on variolation. Variolation was controversial because it usually resulted in illness that required the treated person be quarantined for a few weeks, making it impractical for all but the wealthy, and there was always a risk that procedure could result in death. Despite these drawbacks, there was campaign to promote variolation in the name of progress (Daston, 1998, pp. 83-89) and Daniel Bernoulli developed a mathematical model of the benefits of the technique in 1759, which he presented to the French Academy of Sciences in 1760.

D'Alembert criticised Bernoulli's model on the basis that, in working with expectations across a population, the particular circumstances of an individual were ignored. D'Alembert questioned whether the one in two hundred risk of dying young from the variolation was worth four more years of life expectancy to an individual, though he accepted there was a social benefit in variolation (Daston, 1979, pp. 272-273) (Colombo & Diamanti, 2015, pp. 188-189). D'Alembert developed his criticism by observing that people perceived small clear and present dangers as more threatening than substantial vague and future risks and that lifespans should be measured not just in terms of physical span but also "real" lifespans that measured the time people lived fully and usefully.

While few rejected Bernoulli's model based on d'Alembert's criticism, contemporaries, in particular Condorcet and Laplace, acknowledged that d'Alembert's argument had merit. This was reflected in the description of probability given in the 1765 edition of the *Encyclopédie* (Lubières, 2008). There, following Jacob Bernoulli, a physical approach to probability is described based on the 'nature of things' while the practical approach to probability is founded on experience on the past which is used to predict the future. However, after Poisson most commentators shared Todhunter's 1865 opinion that regarded d'Alembert's position, rooted in a classical approach to uncertainty, as "absurd" to some degree (Daston, 1979, pp. 274-276) and the 'moral' approach to judgement would be replaced by the consequentialist morality of seeking the 'greatest good to the greatest number' founded on universal, mathematical, laws.

## Conclusion

Calvin insisted on an omniscient and omnipotent God but, as the Epicureans had pointed out, individuals' experience of suffering undermined this premise. Calvin's solution was to separate subjective predestination from objective Providence and to employ the Epicurean idea that God's will was hidden. This resulted in a relationship between Divine Providence and chance summarised in the aphorism

> That nothing happens by chance, though the causes may be concealed, but by the will of God; by his secret will which we are unable to explore, but adore with reverence, and by his will which is conveyed to us in the Law and in the Gospel. (Calvin, 1846, Aphorism 13)

Calvin begins by asserting that events are not random, but unlike Stoic or Scholastic philosophers, he also emphasises that God's will was not discernible meaning that mundane experience resulted in the appearance of chance. This was revolutionary in that earlier approaches to *tyche/fortuna* were that it was an irrational motivating power or that chance and suffering were a consequence of a lack of Faith guiding human reason. The idea of God's hidden will appears to be fundamental to Pascal and Fermat's, and Huygens, solution to the *Problem of Points*, noting that the solution relies on considering parts of the path towards a terminal destination that where not part of mundane experience but 'hidden'.

Calvin's doctrine reorientated attitudes to knowledge by removing the distinctions between *phronesis/prudentia* and *episteme/scientia*. Classical approaches to *tyche/fortuna*, which persisted into the eighteenth century, regarded chance events as being highly particular and had to be managed through individual *phronesis*/*prudentia*, personal *elpis*/*spes* or by relying on *pistis*/*fides*. Catholic doctrine was that Grace guided Christian *prudentia* to the universal aim of salvation but could be diverted by individual free will. Because mathematics is concerned with generalisations, neither Pagans nor Christians before Calvin considered universal mathematics relevant to managing uncertain events. Calvin's unification of different types of knowledge into knowledge of God and personal experience enabled the application of universal mathematics to the particularities of judgement in an uncertain world. The change initiated by Calvin is evident in the difference between Machiavelli (Machiavelli, 1909, XVIII) (Randall, 2011, pp. 209-210), who focussed on individual *virtu*/*phronesis,* and Hobbes, who rejected the idea that personal *prudentia* could be part of philosophy because it was individual and not universal (Hobbes, 2017, XLVI). Later, Jacob Bernoulli sub-categorised subjective probability as part of mathematics. The diminution of *phronesis*/*prudentia* was recognised in Daston (1980) (1998) but the significance of Calvin in the process was not recognised.

The initial motivation for this paper was to find out why Huygens considered using the Latin for Hope, *spes*, to denote mathematical expectation, while it is common for the French employ *espérance* today. An explanation for this observation is in how Calvinist theology changed Christian Hope from the Scholastic "uncertain certainty" into a "certain expectation of future happiness" that represented a universal, objective and indubitable idea, not the subjective hope of the pagan or Catholic Hope in salvation, which was contingent on individual will (Aquinas, 2017, II-II.18.4). Being a universal, objective, and indubitable idea, Huygens and Cramer would have recognised *spes*/*espérance* as having the characteristics necessary of a mathematical concept. The presence of *spes*/*espérance* in mathematics is obvious evidence of Calvin's influence on mathematics, as there is no other convincing explanation.

Kendall (1956, pp. 11-12) remarks on Stoic determinism and the role of Providence in Reformed doctrine but does not recognise the Epicurean aspect of Calvin. Puritanism, focussing on Providentialism, re-asserted Stoic ideas and but paid attention to 'signs as evidence'. Spinoza, unconstrained by Calvinist doctrine, asserted that it was possible for people to come to comprehend God's hidden reason, contradicting Calvin. These, more Stoic approaches, set the basis for Hacking's comments that probability emerged just as ideas of divine foreknowledge were being reinforced (2006, pp. 2-3). Hacking also suggests the idea of 'sign as evidence' comes from medical practice. However, the transfer of ideas from medicine to the pioneers of mathematical probability is tenuous. Cardano and Arbuthnot are the only two who had direct medical training and neither of their contributions are as significant as those of Pascal,

Fermat or Huygens, all strongly influenced by Calvinism. Franklin's (2001, pp. 348-352) assertion that the use of *spes* in relation to probability was motivated by the legal term *spei* seems built on a homophone rather than the meaning of the words.

While the importance of Calvin is the central result of this paper, a broader implication is the significance of Epicurean philosophy in the development of mathematical probability. Epicureans are generally regarded as antipathetic to mathematics (Netz, 2015), and so this observation is remarkable. Placing probability theory on a foundation that involves Augustinian, Stoic and Epicurean philosophy, which this paper does, embeds the mathematical theory into a much broader context than just the development of probability in relation to commerce and natural sciences. The significance of classical philosophy suggests Gay's (1966-1969) thesis, which is receiving more interest forty years after its articulation (Edelstein, 2012) (Dew, 2013).

The paper also shows that this 'Epicurean window' was narrow, with the Calvinism of the Dutch and the French (Huguenot) Reformed Churches being overshadowed by British Providential Puritanism and a resurgence of Stoicism, exhibited in Spinoza, from the late seventeenth century. This observation raises the question whether anyone other than Pascal and Huygens, raised under the influence of the Canons of Dort, could have initiated the idea of mathematical expectation as a universal approach to uncertainty, noting that Jacob Bernoulli appears to have held more conventional views.

In terms of contemporary mathematics, this paper argues that the mathematisation of chance is founded on a change in attitude that involved dismissing the importance of individual *phronesis*/*prudentia*. We argue that Calvin initiated this change, but we also observe that the idea was taken up, beyond Calvinism, in the seventeenth century, notably with Hobbes. We also see that d'Alembert resisted the mathematisation of judgements in his opposition to smallpox variolation, since it ignored individual circumstances. Today, when there is public scepticism towards policy founded on mathematics, such as Covid pandemic policies, it is worth bearing in mind that mathematics, necessarily, ignores individual circumstances and so individual, phronetic, attitudes might be passionately opposed to mathematically based, universal, policies. It could benefit mathematicians and statisticians to be conscious of this cost of the mathematisation of judgements.

Simultaneously, the application of mathematics to individual's lives in the form of machine learning is personalising model responses. This is apparent in insurance where the idea of pooling, with the problem of adverse selection, is being replaced by individual risk predictions and this is raising ethical questions about privacy and access to insurance. Current application of mathematics is enacting the ethical problem Calvin grappled with in 1532, of a fully determined society where humans have no effective role.

Being able to illuminate the symmetric problems of depersonalised policies and individual determinism, two significant issues facing contemporary mathematics that are apparently contradictory, on the basis of Calvin's theology is strong evidence of the significance Calvin to mathematics.

## Bibliography

Aeschylus. (1926). Prometheus Bound. In H. W. Smyth (Ed.), *Loeb Classical Library Volumes 145 & 146* (H. W. Smyth, Trans.). Harvard Universrity Press. Retrieved August 2019, from https://www.theoi.com/Text/AeschylusPrometheus.html

Ahnert, T. (2015). *The Moral Culture of the Scottish Enlightenment: 1690-1805.* Yale University Press.

Annas, J. (1995). Prudence and Morality in Ancient and Modern Ethics. *Ethics, 105*(2), 241-257.

Aquinas, T. (2005). *Of God and His Creatures: An Annotated Translation of The Summa Contra Gentiles of St Thomas Aquinas.* (J. Rickaby, Trans.) The Catholic Primer. Retrieved December 2019, from http://anucs.weblogs.anu.edu.au/files/2013/11/St.-Thomas-Aquinas-The-Summa-Contra-Gentiles.pdf

Aquinas, T. (2017). *The Summa Theologiæ of St. Thomas Aquinas.* (K. Knight, Ed., & F. o. Province, Trans.) New Advent. Retrieved August 2019, from New Advent: www.newadvent.org/summa/

Aristotle. (1926). *The Art of Rhetoric.* (J. H. Freese, Ed., & J. H. Freese, Trans.) W. Heinemann. Retrieved September 2019, from data.perseus.org/texts/urn:cts:greekLit:tlg0086.tlg038.perseus-eng1

Aristotle. (2011). *Nicomachean Ethics: Translation, Introduction, Commentary.* (S. Broadie, & C. Rowe, Eds.) Oxford University Press.

Arnauld, A., & Nicole, P. (1850). *Logic, or, The Art of Thinking: being the Port-Royal Logic.* (T. Baynes, Trans.) Sutherland and Knox. Retrieved May 2020, from https://archive.org/details/artofthinking00arnauoft

Assembly of Divines. (2005). *The Westminster Confession of Faith and Longer and Shorter Catechisms.* Orthodox Presbyterian Church. Retrieved June 2020, from https://opc.org/confessions.html

Augustine. (2000). *The City of God Against the Pagans.* (M. Dods, Trans.) Modern Library Classics.

Barner, K. (2001, December). Pierre de Fermat: His life besides mathematics. *Newsletter of the European Mathematical Society*(42).

Barth, K. (1960). *Church Dogmatics.* (G. W. Bromiley, & T. F. Torrance, Trans.) T.& T.Clark.

Bellhouse, D. R. (1988). Probability in the Sixteenth and Seventeenth Centuries: An Analysis of Puritan Casuistry. *International Statistical Review / International Journal of Statistics, 56*(1), 63-74.

Bellhouse, D. R. (1989). A Manuscript on Chance Written by John Arbuthnot. *International Statistical Review, 57*(3), 249-259.

Bellhouse, D. R. (2005). Decoding Cardano's Liber de Ludo Aleae. *Historia Mathematica, 32*, 180-202.

Bernoulli, J. (2005). The Art of Conjecturing, Part 4. In O. Shenyin, *On the Law of Large Numbers* (O. Shenyin, Trans.). Retrieved September 2019, from www.sheynin.de/download/bernoulli.pdf

Bracali, M. (2008). Introduction. In G. Cardano, *De Sapientia libri quinque.* L.S. Olschki.

Bromiley, G. W. (1953). The Doctrine of Christian Hope in Calvin's Institutes. *The Churchman, 67*(3), 148-152.

Brown, R. (1987). History versus Hacking on probability. *History of European Ideas, 8*(6), 655-672.

Brunner, E. (1949). *The Christian Doctrine of God: Dogmatics, Volume I.* Lutherworth Press.

Calvin, J. (1536). *Christianae religionis institutio.* Basel. Retrieved September 2019, from www.e-rara.ch/bau_1/doi/10.3931/e-rara-7379

Calvin, J. (1846). *The Institutes of the Christian Religion.* (H. Beveridge, Trans.) Calvin. Retrieved July 2019, from oll.libertyfund.org/titles/calvin-the-institutes-of-the-christian-religion

Calvin, J. (1859). *Institution de la religion chrestienne.* C. Meyrueis.

Calvin, J. (1969). *Calvin's commentary on Seneca's De Clementia.* (F. Battles, & A. Hugo, Trans.) E.J. Brill. Retrieved August 2019, from media.sabda.org/alkitab-7/LIBRARY/CALVIN/CAL_SENE.PDF

Cambers, A. (2007). Reading, the Godly, and Self-Writing in England, circa 1580–1720. *Jounal of British Studies, 46*(4), 796-825.

Canziani, G. (1992). «Sapientia» e «Prudentia» nella Filosofia Morale di Cardano. *Rivista di Storia della Filosofia, 47*(2), 295-335.

Cape, R. (2003). Cicero and the Development of Prudential Practice at Rome. In R. Hariman (Ed.), *Prudence: classical Virtue, Post-Modern Practice* (pp. 35-65). Pennsylvania State University Press.

Cardano, G. (1953). Liber Ludo Alae. In O. Ore, *Cardano, the Gambling Scholar.* Dover.

Chabbert, P. (1967). Fermat à Castres. *Revue d'histoire des sciences, 20*(4), 337-348.

Clark, M. E. (1983). Spes in the Early Imperial Cult: "The Hope of Augustus". *Numen, 30*(1), 80-105.

Colombo, C., & Diamanti, M. (2015). The smallpox vaccine: the dispute between Bernoulli and d'Alembert and the calculus of probabilities. *Lettera Matematica, 2*(4), 185–192.

Connolly, P. (2014). Locke and Wilkins on Inner Sense and Volition. *Locke Studies, 14*, 239-259.

Cox, V. (1997). Machiavelli and the Rhetorica ad Herennium: Deliberative Rhetoric in The Prince. *The Sixteenth Century Journal, 28*(4), 1109-1141.

Crosby, A. W. (1997). *The Measure of Reality.* Cambridge University Press.

d'Alembert, J.-B. (2009). Preliminary Discourse. In *The Encyclopedia of Diderot & d'Alembert* (R. N. Schwab, & W. E. Rex, Trans.). Michigan Publishing. Retrieved September 2019, from hdl.handle.net/2027/spo.did2222.0001.083

Daston, L. J. (1979). D'Alembert's critique of probability theory. *Historia Mathematica, 6*, 259-279.

Daston, L. J. (1980). Probabilistic Expectation and Rationality in Classical Probability Theory. *Historia Mathematica*, 234-260.

Daston, L. J. (1998). *Classical Probability in the Enlightenment.* Princeton University Press.

David, F. N. (1998). *Games, Gods and Gambling, A history of Probability and Statistical Ideas.* Dover.

de Petris, P. (2008). *Calvin's "Theodicy" in his "Sermons on Job" and the hiddenness of God.* McGill University.

DeBrabander, F. (2007). *Spinoza and the Stoics: Power, Politics and the Passions.* Continuum.

den Uyl, D. J. (1991). *The Virtue of Prudence.* Peter Lang.

Descartes, R. (2008). *A Discourse on Method.* (I. Newby, & G. Newby, Eds.) Project Gutenburg.

Dew, B. (2013). Epicurean and Stoic Enlightenments: The Return of Modern Paganism? *History Compass, 11*(6), 486-495.

Di Somma, E. (2018). *Fides and Secularity: Beyond Charles Taylor's Open Faith.* Wipf and Stock.

Dunlop, A. I. (1992). Provision for Ministers' Widows in Scotland - Eighteenth Century. In A. I. Dunlop (Ed.), *The Scottish Ministers' Widows' Fund 1743-1993.* St. Andrews Press.

Ebbesen, S. (2004). Where Were the Stoics in the Late Middle Ages? In S. K. Stange, & J. Zupko (Eds.), *Stoicism: Traditions and Transformations* (pp. 108-131). Canbridge University Press.

Edelstein, D. (2012). The Classical Turn in Enlightenment Studies. *Modern Intellectual History, 9*(1), 61-71.

Essary, K. (2017). *Erasmus and Calvin on the foolishness of God: Reason and Emotion in the Christian Philosophy.* University of Toronto Press.

Euripides. (2010). *Suppliant Women Ικέτιδες.* (G. Theodoridis, Trans.) Retrieved September 2019, from Bacchicstage: bacchicstage.wordpress.com/euripides/suppliant-women/

Faraguna, M. (2012). Pistis and apistia: aspects of the development of social and economic relations in classical Greece. *Mediterraneo antico, 15*(1-2), 355-374.

Franklin, J. (2001). *The Science of Conjecture: Evidence and Probability before Pascal.* Johns Hopkins University Press.

Franklin, J. (2006). *Catholic Values and Australian Realities.* Connor Court Publishing .

Franklin, J. (2009). Science by Conceptual Analysis: The Genius of the Late Scholastics. *Studia Neoaristotelica, 6*, 209–233.

Garber, D., & Zabell, S. (1979). On the Emergence of Probability. *Archive for History of Exact Sciences, 21*(1), 33-53.

Gay, P. (1966-1969). *The Enlightenment: An Interpretation.* Alfred A. Knopf.

Gilby, E. (2009). The Language of Fortune in Descartes. In J. D. Lyons, & K. Wine (Ed.), *Chance, Literature, and Culture in Early Modern France.* Ashgate.

Goetzman, W. (2016). *Money Changes Everything: How Finance Made Civilisation Possible.* Princeton University Press.

Gootje, N. H. (2005). Calvin on Epicurus and the Epicureans. *Calvin Theological Journal, 40*, 33-48.

Hacking, I. (2006). *The emergence of probability.* Cambridge University Press.
35

Hadden, R. W. (1994). *On the Shoulders of Merchants: Exchange and the Mathematical Conception of Nature in Early Modern Europe.* State University of New York Press.

Hald, A. (1990). *A History of Probability and Statistics and their Applications before 1750.* Wiley.

Hare, D. J., & Scott, W. F. (1992). The Scottish Ministers' Widows' Fund of 1744. In A. I. Dunop (Ed.), *The Scottish Ministers' Widows' Fund 1743--1993.* St. Andrews Press.

Hariman, R. (2003). Preface. In R. Hariman (Ed.), *Prudence: Classical Virtue, Postmodern Practice* (pp. vii-x). Pennsylvania University Press.

Hay, D. M. (1989). Pistis as "Ground for Faith" in Hellenized Judaism and Paul. *Journal of Biblical Literature, 108*(3), 461-476.

Heims, S. J. (1980). *John von Neumann and Norbert Weiner: From Mathematicians to the Technologies of Life and Death.* MIT Press.

Hobbes, T. (2017). *Leviathan.* (J. Bennett, Ed.) Early Modern Texts. Retrieved 2019 September, from www.earlymoderntexts.com/authors/hobbes

Howson, C. (1978). The Prehistory of Chance. *The British Journal for the Philosophy of Science, 29*(3), 274-280.

Huygens, C. (1714). *The Value of All Chances in games of Fortune.* T. Woodward. Retrieved September 2019, from www.york.ac.uk/depts/maths/histstat/huygens.pdf

Ingham, M. B. (2016). Stoic Influences in the Later Middle Ages . In J. Sellers (Ed.), *The Routledge Handbook of the Stoic Tradition* (pp. 120-136). Routledge.

Kaye, J. (1998). *Economy and Nature in the Fourteenth Century.* Cambridge University Press.

Kemp, C. (2014). The Real 'Letter to Arbuthnot'? a Motive For Hume's Probability Theory in an Early Modern Design Argument. *British Journal for the History of Philosophy, 22*(3), 468-491.

Kendall, M. G. (1956). Studies in the History of Probability and Statistics: II. The Beginnings of a Probability Calculus. *Biometrica, 43*(1/2), 1-14.

Kershaw, S. P. (1986). *Personification in the Hellenistic world : Tyche, Kairos, Nemesis.* Thesis, University of Bristol.

Lactantius. (n.d.). *On the Anger of God.* Retrieved September 2019, from www.newadvent.org/fathers/0703.htm

Le Goff, J. (1990). *The Birth of Purgatory.* (A. Goldhammer, Trans.) Scolar Press.

Lee, P. M. (n.d.). *Fermat and Pascal on Probability.* University of York. Retrieved September 2019, from www.york.ac.uk/depts/maths/histstat/pascal.pdf

Legros, A. (2009). Montaigne Between Fortune and Providence. In J. D. Lyons, & K. Wine (Ed.), *Chance, Literature, and Culture in Early Modern France* (K. F. McConnell, Trans.). Ashgate.

Lewin, C. ( 2004). Graunt, John (1620–1674), statistician. *Oxford Dictionary of National Biography*.

Locke, J. (2017b). *An Essay Concerning Human Understanding, Book IV.* (J. Bennett, Ed.) Retrieved September 2019, from Early Modern Texts: www.earlymoderntexts.com/assets/pdfs/locke1690book4.pdf

Long, A. A. (2006). Chance and Laws of Nature in Epicureanism. In A. Long, *From Epicurus to Epictetus: Studies in Hellenistic and Roman Philosophy* (pp. 157-177). Oxford University Press.

Lubières, C.-B. (2008). Probability. In Diderot, & d'Alembert, *The Encyclopedia Collaborative Translation Project* (D. Weiner., Trans.). University of Michigan Library. Retrieved January 14, 2019, from http://hdl.handle.net/2027/spo.did2222.0000.983

Luehrs, R. (1987). Population and Utopia in the Thought of Robert Wallace. *Eighteenth-Century Studies, 20*(3), 313-335.

MacDonald, P. S. (2002). Descartes: The Lost Episodes. *Journal of the History of Philosophy, 40*(4).

Machiavelli, N. (1909). *The Prince.* (L. Ricci, Trans.) Oxford University Press. Retrieved from www.gutenberg.org/ebooks/57037

McKeveley, C. (2018). *Reading and Writing Certainty in Early Modern England.* Southern Methodist University.

Meyns, C. (2017). Leibniz and Probability in the Moral Domain. In S. L, V. E., & W. J. (Eds.), *Tercentenary Essays on the Philosophy and Science of Leibniz.* Palgrave Macmillan.

Miller, J. (2015). *Spinoza and the Stoics.* Cambridge University Press.

Momigliano, A. (1984). Religion in Athens Rome and Jerusalem in the First Century B. C. *Annali della Scuola Normale Superiore di Pisa. Classe di Lettere e Filosofia, Serie III, 14*(3), 873-892.

Monsalve, F. (2014). Scholastic just price versus current market price: is it merely a matter of labelling? *The European Journal of the History of Economic Thought, 21*(1), 4-20.

Montaigne, M. (1877). *The Essays of Montaigne.* (W. C. Hazlitt, Ed., & C. Cotton, Trans.) Reeves and Turner. Retrieved September 2019, from www.gutenberg.org/files/3600/3600-0.txt

Morgan, T. (2015). *Roman Faith and Christian Faith: Pistis and Fides in the Early Roman Empire and Early Churches.* Oxford University Press.

Muller, R. (1986). *Christ and the Decree: Christology and Predestination in Reformed Theology from Calvin to Perkins.* The Labyrinth Press.

Netz, R. (2015). Were There Epicurean Mathematicians? In B. Inwood (Ed.), *Oxford Studies in Ancient Philosophy* (Vol. 49). Oxford University Press.

Nussbaum, M. (2001). *The Fragility of Goodness: Luck and Ethics in Greek Tragedy and Philosophy.* Cambridge University Press.

Ong, M.-C. (2003). *John Calvin on providence: the locus classicus in context.* King's College London.

Oxenham, S. (2000). *A touchstone the written word' : experimental Calvanist life-writing and the anxiety of.* King's College London.

Partree, J. (2005). *Calvin and Classical Philosophy.* Westminster John Knox Press.

Pascal, B. (1958). *Pensées.* E. P. Dutton & Co. Retrieved December 2019, from https://www.gutenberg.org/files/18269/18269-h/18269-h.htm

Peterson, D. J. (1994). *Political economy in transition: From classical humanism to commercial society - Robert Wallace of Edinburgh.* PhD Thesis, University of Illinois at Urbana-Champaign.

Pitkin, B. (2016). Erasmus, Calvin, and the Faces of Stoicism in Renaissance and Reformation Thought. In J. Sellers (Ed.), *The Routledge Handbook of the Stoic Tradition* (pp. 168-183). Routledge.

Randall, D. (2011). The Prudential Public Sphere. *Philosophy & Rhetoric, 44*(3), 205-226.

Reardon, P. H. (1975). Calvin on Providence: The Developement of an Insight. *Scottish Journal of Theology, 28*, 517-534.

Reeves, J. (2015). The Secularization of Chance: Toward Understanding the Impact of the Probability Revolution on Christian Belief in Divine Providence. *Zygon, 50*(3), 604-620.

Regosin, R. (2009). Prudence and the Ethics of Contingency in Montaigne's Essais. In J. D. Lyons, & K. Wine (Ed.), *Chance, Literature, and Culture in Early Modern France.* Ashgate.

Restivo, S. (1982). Mathematics and the Sociology of Knowledge. *Science Communication, 4*(1), 127-144.

Ross, A. (2004). Arbuthnot [Arbuthnott], John (bap. 1667, d. 1735), physician and satirist. *Oxford Dictionary of National Biography*.

Ryrie, A. (2016). The nature of spiritual experience. In U. Rublack (Ed.), *The Oxford handbook of the Protestant Reformations* (pp. 47-63). Oxford University Press.

Shoesmith, E. (1987). The Continental Controversy over Arbuthnot's Argument for Divine Providence. *Historia Mathematica, 14*, 133–46.

Shulman, G. (1988). Hobbes, Puritans, and Promethean Politics. *Political Theory, 16*(3), 426-443.

Spinoza, B. (1901). Ethics. George Bell and Sons. Retrieved September 2019, from www.gutenberg.org/files/3800/3800-h/3800-h.htm

Stachniewski, J. (2008). Introduction. In J. Bunyan, J. Stachniewski, & A. Pacheco (Eds.), *Grace Abounding with Other Spiritual Autobiographies.* Oxford University Press.

Starr, G. A. (1965). *Defoe and Spiritual Autobiography.* Princeton University Press.

Steinmetz, D. (2010). *Calvin in Context.* Oxford University Press.

Stoffele, B. (2006). *Christiaan Huygens – A family affair: Fashioning a family in early-modern court-culture.* MSc thesis, Utrecht University, History and Philosophy of Science.

Strohm, C. (2019). Sixteenth-Century French Legal Education and Calvin's Legal Education. In B. C. Brewer, & D. M. Whitford (Eds.), *Calvin and the Early Reformation* (pp. 44-57). Brill.

Sylla, E. D. (2006). Commercial arithmetic, theology and the intellectual foundations of Jacob Bernoulli's Art of Conjecturing. In G. Poitras (Ed.), *Pioneers of Financial Economics: contributions prior to Irvin Fisher* (pp. 11-45). Edward Elgar.

Synod of Dort. (1619). The Canons of Dort. *The National Synod of the National Reformed Chuch.* Dordrecht. Retrieved December 2019, from https://prts.edu/wp-content/uploads/2016/12/Canons-of-Dort-with-Intro.pdf

Theognis. (1931). *Elegy and Iambus.* (J. Edmonds, Trans.) Harvard University Press. Retrieved August 2019, from demonax.info/doku.php?id=text:theognis_poems

Todhunter, I. (2014). *Todhunter. (2014). A history of the mathematical theory of probability : from the time of Pascal to that of Laplace .* Cambridge University Press.

Torrance Kirby, W. J. (2003). Stoic and Epicurean? Calvin's Dialectical Account of Providence in the Institute. *International Journal of Systematic Theology, 5*, 309-322.

van Brankel, J. (1976). Some Remarks on the Prehistory of the Concept of Statistical Probability. *Archive for History of Exact Sciences, 16*(2), 119-136.

Verbeke, G. (1983). *The Presence of Stoicism in Medeival Thought.* Catholic University of America Press.

Walsh, P. G. (1974). Spes Romana, Spes Christiana. *Prudentia, 6*, 33-42.

Wilkins, J. (1649). *A discourse concerning the beauty of providence in all the rugged passages of it very seasonable to quiet and support the heart in these times of publick confusion* (Early English books Text Creation Partnership ed.). University of Michigan. Retrieved June 2020, from https://quod.lib.umich.edu/e/eebo/A66025.0001.001

Zimmermann, R. (1996). *The Law of Obligations: Roman Foundations of the Civilian Tradition.* Oxford University Press.